\documentclass[12pt]{article}
\def\C{\mathbf{C}}
\def\bC{\mathbf{\overline{C}}}
\def\N{\mathbf{N}}
\begin{document}
\title{The work of Walter Bergweiler in value distribution of meromorphic
functions}
\author{Alexandre Eremenko}
\maketitle
\begin{abstract}
This is the colloquium talk in CAU Kiel delivered on June 7, 2024 on
the occasion of Walter Bergweiler's retirement.
It contains a survey  of his work in the value distribution
theory of meromorphic functions of one complex variable. The talk was aimed
at a broad audience, so a brief introduction to Nevanlinna theory is included.

\end{abstract}

Value distribution theory begins with the celebrated theorem of Emile Picard:
\vspace{.1in}

\noindent
{\em If an entire function omits two values, then it is constant}.
\vspace{.1in}

This theorem (and its two-lines proof) published in 1879
made an enormous impression, and
they say that all seminar talks in Paris in the early 20th century began with 
the words: ``Picard proved...''.

One can restate it as follows {\em if a meromorphic function in the plane
omits three values then it is constant}.
The generalization known as the Great Picard's theorem says that
{\em If a meromorphic function in a punctured neighborhood of a point
omits three values in $\bC$ then this point is a pole or removable
singularity.}

Schottky \cite{Schottky} proved the
corresponding result for the unit disk. We state it in the form
proposed by Montel: {\em The family of all
meromorphic functions in the unit disk
which omit three values, say $0,1,\infty$ is normal.}
I recall that a family $F$ of meromorphic functions in a region $D$
is called {\em normal}
if from any sequence $(f_n)$ in $F$ one can select a subsequence
which converges uniformly (with respect to the spherical metric in
the image) on every compact subset of $D$.

One can show that the theorems of Picard and Montel can be simply derived
from each other, and this is the basis
of the famous
``Bloch's Principle'': {\em If some
conditions imply that a meromorphic function in the plane is constant,
then the same conditions imposed on a family of meromorphic functions
in any region imply that this family is normal, and conversely.}

This heuristic
principle turned out to be very fruitful, see, for example
Walter's survey \cite{BB}. It is equally important in the theory
of holomorphic curves, \cite{Mc}, \cite{Yam}.

The condition that some value $a$ is omitted can be usually replaced 
by the condition that all solutions of $f(z)=a$ have high multiplicity.
Nevanlinna proved the following generalization of Picard's theorem:
{\em If $f$ is a non-constant
meromorphic function in $\C$,
and $a_j$ have the property that all solutions of $f(z)=a_j$ have
multiplicities
at least $m_j$, then
\begin{equation}
\label{1}
\sum_j\left(1-\frac{1}{m_j}\right)\leq 2,
\end{equation}
where $m_j=\infty$ if $a_j$ is omitted.
}
This result was first obtained by Nevanlinna, as a consequence of
his deep {\em Nevanlinna theory} which can be considered as a far-reaching
quantification of Picard's theorems.
Since its inception it became the main part of the theory of meromorphic
functions, and I will briefly introduce the main notions and results.
Let $f:\C\to\bC$ be a meromorphic function, and
$$A(r,f)=\frac{1}{\pi}\int\int_{x^2+y^2\leq r^2}
\frac{|f'(x+iy)|^2}{(1+|f(x+iy)|^2)^2}dxdy,$$
the area of the image of the disk under $f$ with respect to the spherical
metric, divided by the area of the sphere. So for a rational function
of degree $d$ we have $A(r,f)\to d, r\to\infty$. Then the Nevanlinna
characteristic is defined as an average of $A(r,f)$ with respect to $r$:
$$T(r,f)=\int_{0}^r \frac{A(t,f)}{t}dt.$$
This function measures the growth of a meromorphic function. For a rational
function $T(r,f)\sim d\log r$, and for transcendental meromorphic functions
we have $T(r,f)/\log r\to\infty$. Nevanlinna characteristic has formal
algebraic
properties similar to the degree of a rational function, and this is one
reason why we averaged with respect to $r$. For example,
$$T(r,f)=T(r,1/f),\quad T(r,f^n)=nT(r,f),\quad T(r,fg)\leq T(r,f)+T(r,g)+O(1),$$
and $T(r,f+g)\leq T(r,f)+T(r,g)+O(1)$.

The order of a meromorphic function is defined as
$$\rho:=\limsup_{r\to\infty}\frac{\log T(r,f)}{\log r}\in[0,\infty].$$
So for example, $e^z$ and $\sin z$ are of order $1$, $e^{z^p}$
is of order $p$, $\wp(z)$ has order $2$, and $e^{e^z}$ is of infinite order.

For every $a\in\bC$ we define $n(r,a,f)$ as the number (counting multiplicity)
of solutions of $f(z)=a$ in the disk $|z|\leq r$, and then apply the averaging
in $r$ as above:
$$N(r,a,f)=\int_0^r(n(t,a,f)-n(0,a,f))\frac{dt}{t}.$$
Then the First Main Theorem of Nevanlinna says that for all $a\in\bC$
we have
$$N(r,a,f)\leq T(r,f)+O(1),\quad r\to\infty,$$
and the Second Main Theorem says that for any finite set
$\{a_1,\ldots,a_q\}$
$$\sum_{j=1}^q N(r,a_j,f)\geq (q-2+o(1))T(r,f)+N_1(r,f),$$
when $r\to\infty$ avoiding a set of finite measure.
Here $N_1(r,f)$ is the ramification term: $n_1(r,f)$ is the number
of critical points (counting multiplicity) in the disk $|z|\leq r$,
and $N_1$ is obtained from $n_1$ by the averaging as above. The meaning of
the Second Main Theorem can be made clearer if we define the deficiency
$$\delta(a,f)=1-\limsup_{r\to\infty}\frac{N(r,a,f)}{T(r,f)}.$$
So $0\leq\delta(a,f)\leq 1$, and if the value $a$
is omitted then $\delta(a,f)=1$. Then the Second Main Theorem implies that
$\delta(a,f)>0$ for at most countable set of $a$ (they are called deficient),
and
\begin{equation}\label{def}
\sum_{a\in\bC}\delta(a,f)\leq 2.
\end{equation}
This is a generalization of Picard's theorem. If one takes multiplicities
into account, then (\ref{1}) follows.

Returning to Picard's theorem, it inspired a lot
of variations and generalizations in 20th century. Our main subject
here will be the analogs of Picard's theorem involving derivatives.

Among the earliest results on this subject,
I mention a theorem of
Milloux \cite{Mi} {\em If an entire function omits some finite value,
and its derivative omits some non-zero value, then $f$ is constant}.
The corresponding normality criterion was proved by Miranda \cite{Miranda}
5 years earlier. Example $e^z$ shows that $f'$ can indeed omit zero.
This shows that the value $0$ is special for the derivative.

On the other hand, Csillag \cite{C} proved that {\em If $f$ is entire and
$f,f^{(m)},f^{(n)}$ omit zero for some $0<m<n$, then $f(z)=e^{az+b}$.}
The case $(m,n)=(1,2)$
was earlier proved by Saxer.

A great breakthrough in these questions was made by Hayman \cite{H}. The main
result of his article is the so-called 
\vspace{.1in}

\noindent
{\bf Hayman's Alternative.} {\em If $f$ is a transcendental meromorphic function
in the plane and $n$ is a positive integer,
then either $f$ assumes every finite value infinitely many times,
or $f^{(n)}$ assumes every finite non-zero value infinitely many times.}

\vspace{.1in}

It is noteworthy that this theorem gives only two conditions
on omitted values which imply that the function it is constant.
Picard, Schottky, Milloux and Csillag theorems when stated
for meromorphic functions
require 3 conditions (they all assume that $\infty$ is omitted).

Hayman's result generalizes the theorem of Milloux to meromorphic functions.

In this paper Hayman proved many other results, but he assumed that they are
not best possible, and stated several conjectures.
The most outstanding conjectures are these two:
\vspace{.1in}

\noindent
{\bf Hayman's Conjecture 1.} {\em If $f$ is a transcendental
meromorphic function in $\C$,
and $n\geq 2$, then $(f^n)'$ assumes all non-zero complex values 
infinitely many times.}
\vspace{.1in}

He himself proved this for $n\geq 4$, and later Mues extended this to $n=3$,
and Clunie proved this conjecture for entire functions.
\vspace{.1in}

\noindent
{\bf Hayman's Conjecture 2.} {\em If $f$ is a meromorphic function in $\C$
and both $f$ and $f^{(n)}$ are free of zeros for some $n\geq 2$,
then $f(z)=e^{az+b}$ or
$f(z)=(az+b)^{-m}$ with som integer $m$.}
\vspace{.1in}

This was proved by Frank for $n\geq 3$ and by Mues for functions
of finite order.
The final result with $n=2$ is due to Langley \cite{L1} who developed
the methods of Mues and Frank in a very subtle way.
\vspace{.1in}

I met Walter in winter 1992/3 in a conference in Joensuu. For me this was
probably the first international conference that I was able to attend (I moved
from the Soviet Union to the US in 1991). 
We shared a 3-bedroom apartment with my adviser A. A. Goldberg and for him
this was also one of the first trips abroad. We had no common language
for all three of us: Goldberg did not speak English, so he had to
communicate with Walter in German and with me in Russian. And I translated
his conference talk into English.

Next summer I was invited by Walter to stay with him and his wife in Aachen,
where he had his first position. When meeting me in the train station,
he explained to me his proposal to prove Hayman's Conjecture 1 in
the remaining case $n=2$ {\em for functions of finite order}.

The result was the paper \cite{BE} which we wrote that summer and submitted
to the journal. The proof was unusual for the theory of meromorphic functions
since it used holomorphic dynamics.
Holomorphic dynamics which was originated by Fatou (with some result
obtained independently by Julia, Latt\`es and Ritt) in the beginning
of 20th century experienced a strong revival since 1982.

We were both interested in holomorphic dynamics at that time: I was
deeply impressed with the work of Misha Lyubich in the early 1980s,
and we wrote several joint papers in Kharkiv on dynamics of transcendental
functions. Walter started thinking about holomorphic dynamics under the
influence of J. Hubbard, when he visited Cornell University as a postdoc
in 1987-89.
By 1993 we both wrote surveys on holomorphic dynamics: my survey was
joint with Lyubich, and Walter's survey of dynamics of transcendental
meromorphic function is to this time the best introductory survey;
it accumulated 417 citations on MSN by the time I write this.

The idea that Walter explained to me in the train station was the following:
if $f$ has finitely many zeros, then Hayman's alternative is applicable,
and we are done. So assume that it has infinitely many zeros. They are
all multiple zeros of $f^n$. Then the function
$$g(z)=z-f^n(z)$$
has infinitely many fixed points, and these fixed points are neutral,
with multiplier $1$, that is $g'(a)=1$ at all these fixed points $a$.
To each such fixed point an immediate domain
of attraction $P$ is associated (the so called Leau's petal). 
This domain is invariant under $g$, and according to the fundamental theorem
of Fatou, $P$ must contain either critical
of asymptotic value of $g$. (We say that $a\in\bC$ is an asymptotic value
of a function $g$ meromorphic in $\C$ if here is a curve tending
to $\infty$ such that $g(z)\to a$ along this curve.)
If this singular value which exists by Fatou's theorem
is critical then $g'(z)=0$
at some point $z$ in the domain of attraction of $P$, and since these
domains are disjoint for different neutral fixed points,
we obtain infinitely
many critical points, which are solutions of $(f^{n})'=1$.

So it remains to deal with the possibility that $g$ has infinitely
many asymptotic values but only finitely many critical values.
Such meromorphic functions exist, an example can be found in \cite{V}.

But it turns out that such thing is impossible for functions
of {\em finite order}. This was a great insight of Walter which he explained
to me during that conversation in the train station. It took few weeks
to work out the details, and the final result was the following
\vspace{.1in}

\noindent
{\bf Theorem 1.} {\em Let $f$ be a meromorphic function of finite order
$\rho$.
Then every asymptotic value of $f$, except at most $2\rho$ of them,
is a limit point of critical values different from this asymptotic value.}
\vspace{.1in}

By using the dynamics argument outlined above, we obtain
\vspace{.1in}

\noindent
{\bf Theorem 2.} {\em Let $f$ be a transcendental
meromorphic function of finite order,
and assume that it has infinitely many multiple zeros. Then $f'$
assumes every finite non-zero value infinitely often.}
\vspace{.1in}

Unfortunately, this theorem fails for functions of infinite order,
so it seemed that we could not prove Hayman's Conjecture 1 in full generality
with this method.

I knew from the work of Zalcman \cite{Z} and Brody that there exists
a marvelous device to reduce Picard-type theorems to their special cases
for functions of finite order. This is called the Zalcman Rescaling Lemma
to the specialists in classical function theory, and Brody's Rescaling Lemma
to the specialists in holomorphic curves. Actually this argument goes back
to the work of Bloch and Valiron in the beginning of 20th century.
However this rescaling argument did not seem to be appropriate for 
working with derivatives.

Only next summer (1994) when browsing the book by Schiff on normal families,
I learned about the paper by X. Pang \cite{Pang} of 1989 where a suitable
generalization of Zalcman's Lemma was proved\footnote{Independently of Pang,
similar arguments were used in Schwick's thesis \cite{Schwick}.}. I immediately informed
Walter (and Zalcman) and since our paper was not yet published, we managed
to revise it, and obtain the full proof of Hayman's conjecture \cite{BE}.
Our final result was
\vspace{.1in}

\noindent
{\bf Theorem 3}. {\em For every transcendental meromorphic function
$f$, and every positive integers $k<n$, the function $(f^n)^{(k)}$
takes every finite non-zero value infinitely often.}
\vspace{.1in}

Meanwhile, the statement of our Theorem 2 reached China (this was before
the advent of the arXiv): Walter gave a conference
talk in Nanjing in June 1994, and the consequence was an ``independent'' proof
of Hayman's conjecture published by Huaihui Chen and Mingliang Fang
\cite{Chen}. According to a later personal communication
from X. Pang,
he knew how to reduce Hayman's conjecture to its special case of finite order
even before our paper was written, but Pang was not present at Walter's
talk in Nanjing.

In a later paper \cite{BN}, Walter applied Pang's lemma to give a new
simple proof of Langley's theorem stated above. The corresponding normality criterion
is
\vspace{.1in}

{\em Let $F$ be a family of meromorphic functions $f$ such that $ff''$ omit
zero. Then the family $\{ f'/f:f\in F\}$ is normal.}
\vspace{.1in}

Since then our Theorem 3 was very much used and generalized,
paper \cite{BE} has 230 citations on MSN and 550 on Google Scholar.

The ultimate
generalization for the first derivative was conjectured by Walter
and proved by  Jianming Chang \cite{Chang}, who used a deep intermediate
result
of Nevo, Pang and Zalcman \cite{NPZ}:
\vspace{.1in}

\noindent
{\em Suppose that for a transcendental meromorphic function $f$, $|f'|$
is bounded on
the zero set of $f$. Then $f$ takes every finite non-zero value
infinitely often.}
\vspace{.1in}

There is an essential difference between this and Theorem 2:
(which is not true for functions of infinite order):
In Theorem 2 the assumption is that there are infinitely many multiple zeros,
while in the theorem of Chang the assumption of boundedness
of $|f'|$ applies to {\em all} zeros. 

Despite many variations and generalizations, no new proof of Theorem 3
was found until 2020: all variants and generalizations used our Theorem 3.
I always felt that our proof was
unsatisfactory from the philosophical
point of view. One reason for this was the use of dynamics: why a non-dynamical
theorem should use dynamics in its proof?
Walter showed me how dynamics can be
eliminated (see \cite{EB} where his argument is recorded). It uses instead
the so-called Logarithmic Change of the Variable, which is a common tool 
in transcendental dynamics, but formally does not belong to it.
(It was invented by Teichm\"uller for a problem of value distribution
theory).

More important philosophical objection is related to the weird overall
scheme: we prove Theorem 2 for functions of finite order,
({\em and it is not true without this restriction!}). And then by some
magic trick we derive Theorem 3 which applies to all meromorphic functions.

The new proof, which uses neither our arguments nor Pang's lemma,
was obtained recently in \cite{An} as a corollary of the very deep result
of Yamanoi \cite{Y}, and I briefly describe his result.

Erwin Mues conjectured in 1971 that the defect relation for 
the derivative of a meromorphic function can be improved to
$$\sum_{a\neq\infty} \delta(a,f')\leq 1.$$
In his attempts to prove this conjecture Goldberg made a new one
{\em For all meromorphic functions
\begin{equation}\label{G}
\overline{N}(r,\infty,f)\leq N(r,0,f'')+o(T(r,f)),
\end{equation}
(as $r\to\infty$ avoiding a small exceptional set),
where $\overline{N}$ is the averaged counting function of distinct
poles (that is counted without multiplicity).} Goldberg showed that
this implies Mues's conjecture, and proved (\ref{G}) for the case
when all poles of $f$ are simple.
Actually, both the Mues conjecture and Goldberg's
conjecture for the case of simple poles
can be derived from the arguments of Frank and Weissenborn \cite{F}.

It is amazing that no perturbation argument seems to derive the
Mues and Goldberg
conjectures from this special ``generic'' case.

Langley \cite{L4} proved that for meromorphic functions of finite order
the condition that $f^{(k)}$ for some $k\geq 2$
has finitely many zeros implies that $f$
has finitely many poles. However he constructed an example of a 
meromorphic function of arbitrarily slow infinite order of growth
for which $f''$ has no zeros but $f$ has infinitely many poles.
So there is no Picard-type theorem corresponding to Goldberg's
conjecture!

Goldberg's conjecture was proved by Yamanoi \cite{Y} in the following
stronger form
$$k\overline{N}(r,\infty,f)+
\sum_{a\in A}N_1(r,a,f)\leq N(r,0,f^{(k+1)})+o(T(r,f))$$
outside a set of $r$ of zero logarithmic density, where $k$ is
any positive integer and  $A$ is any
finite set.
This
77 pages paper is a real tour de force; this is probably the most
difficult of all results in value distribution of meromorphic functions
in dimension one. In the recent decades, the main interests in
Nevanlinna theory shifted
from meromorphic functions to holomorphic curves, but this result
of Yamanoi (and his other results) on one-dimensional Nevanlinna theory
show that it is still very much alive.

I make few comments on Theorem 1, which is of independent interest.
It is related to the famous Denjoy--Carleman--Ahlfors theorem, and
Iversen's classification of singularities of inverse functions.
If $a$ is an asymptotic value of a meromorphic function $f$, 
and $D(a,r)$ is an arbitrarily small (spherical) disk centered at $a$,
then the preimage $f^{-1}(D(a,r))$ has an unbounded component $V$,
the one that contains the asymptotic curve. If for some $r>0$
the equation $f(z)=a$ has no solutions in $V$, we say that the singularity
of $f^{-1}$ over $a$ is direct. Otherwise it is indirect.
For example $\log z$ has two direct singularities: one over $0$ another
over $\infty$. And $(\sin z)/z$ has two indirect singularities over $0$
and two direct ones over $\infty$.

The famous Denjoy--Carleman--Ahlfors theorem says that for a meromorphic
function of order $\rho$, the inverse has at most $2\rho$ direct singularities.
And the more precise form of our Theorem 1 is:
\vspace{.1in}

{\em If $f$ is a meromorphic function of finite order, then every neighborhood
of an indirect singularity over $a$ (that is every component $V$
in the above explanation) contains infinitely many critical points
with critical values distinct from $a$.}
\vspace{.1in}

Iversen's classification was rarely used in value distribution theory
before our paper, but since its publication it became a common tool in
value distribution and transcendental holomorphic dynamics.

Returning to Walter's work, I would like to
mention here three his other contributions.

In his paper \cite{B1} he proved the following
\vspace{.1in}

\noindent
{\bf Theorem 4.} {\em Let $f$ be a meromorphic function of finite order,
and $a\in\C$. If $a\neq 1$ and $a\neq (n+1)/n$ for all $n\in\N$, 
then $ff''-a(f')^2$ has infinitely many zeros, unless $f=Re^P$ with
rational $R$ and polynomial $P$.}
\vspace{.1in}

He conjectured that the assumption of finite order can be removed in
this theorem; this would give a generalization of Langley's proof
of Hayman's Conjecture 2, which corresponds to $a=0$.
Walter's conjecture was proved by Langley \cite{L2}, without use
of Pang's lemma. A simple derivation based on Pang's Lemma
appeared in \cite{Ln}. 

In the joint paper \cite{BL} Walter and Jim study the question
whether in Hayman's alternative one can replace omitted values by
high multiplicity. As expected, they show that this is so, but 
probably their results are not best possible. Improving a previous result
by Zhen Hua Chen, they show that a meromorphic function $f$ whose
zeros have multiplicity $\geq m$ while $1$-points of $f^{(k)}$
have multiplicity $\geq n$ must be constant if
$$\left(2k+3+\frac{2}{k}\right)/m+\left(2k+2+\frac{2}{k}\right)/n<1.$$
So when $k=1$, we have the condition
$$\frac{7}{m}+\frac{8}{n}<1.$$
They also constructed an example of a non-constant meromorphic function
with $m=2$ and $n=3$.
On the other hand, Nevo, Pang and  Zalcman proved that
\vspace{.1in}

\noindent
The derivative of a transcendental
meromorphic function whose almost all zeros are multiple
takes every non-zero value infinitely often. 
\vspace{.1in}

Notice that the rational function
$$f(z)=\frac{(z-a)^2}{(z-b)}$$
has only one zero in $\C$, this zero is multiple, and $f'$ omits $1$.

So it is an open problem to obtain a best possible
result about multiplicities in Hayman's alternative.
The statement about normality corresponding to this via Bloch's Principle
also holds in this setting.
\vspace{.1in}

In his influential paper \cite{Bloch} Bloch mentioned two philosophical
principles: one that we stated in the beginning and another the
``Principle
of topological continuity''. To explain this second principle, consider a 
disk $V\subset \bC$. A component of preimage $f^{-1}(V)$ under a meromorphic
function $f$ is called an {\em island} over $V$ if this component
is bounded.
In this case the restriction of $f$ on this component is a ramified covering
and thus has a {\em degree} which is a positive integer. This degree
is called the multiplicity of an island. 
If $a\in V$ is a point,
and $f$ has an $a$-point
of multiplicity $m$ in an island over $V$, then the multiplicity of
the island is at least $m$. Bloch proposed the
following generalization of (\ref{1}):

{\em If $V_j$ are
disks with disjoint closures, and the multiplicity of all
islands over $V_j$ of a non-constant  meromorphic function $f$ in $\C$
are at least $m_j$, then (\ref{1}) holds.} 

If there are no islands over $V_j$ we set $m_j=\infty$.

Bloch's conjecture was proved as a consequence of very deep
``\"Uber\-la\-ger\-ungs\-fl\"a\-chen\-the\-orie'' of Ahlfors for which he
obtained one of the first two Fields medals in 1936. This theory is very
technical, and it does not use the Nevanlinna's result (\ref{1}).
A modern proof of (\ref{1}) can be obtained without
Nevanlinna theory, a simple proof was obtained by Robinson.

Several simplified derivations of Ahlfors theory also exist, but none
of them can be called really simple.

Walter found in \cite{BA} a simple derivation of the Ahlfors island
result from (\ref{1}), in complete agreement with Bloch's vision.
So this is by far the simplest available proof of the Ahlfors's islands
theorem.

The ``Principle of Topological Continuity'' as stated by Bloch, is
not a theorem but a philosophical principle, and in our most recent
paper with Walter \cite{BE2} we investigate the cases when
it can be converted into
theorems. For example, one simple consequence of Nevanlinna theory
(proved by H. Selberg) is 
\vspace{.1in}

\noindent
{\em Let $f$ be a non-constant entire function of order less than 1,
and $a,b$ are two points in $\C$.
Then $f$ one of these two points has a simple preimage.}
\vspace{.1in}

Notice that $\cos\sqrt{z}$ is of order $1/2$, and 
has only one simple preimage of the pair of points $\{\pm1\}$.

Bloch stated that to this should correspond:
\vspace{.1in}

\noindent
{\em Let $f$ be a non-constant entire function of order less than $1$,
and $A,B$ are two discs in $\C$ with disjoint closures. Then over one
of these disks there is a simple island.}
\vspace{.1in}

We proved this in \cite{BE2}. But some other similar statements
involving meromorphic functions and simple islands over one of the
4 disks with disjoint closures turn out to be wrong.

\vspace{.1in}

{\em Department of Mathematics

Purdue University

West Lafayette, IN 47907

eremenko@purdue.edu}

\begin{thebibliography}{1}
\bibitem{An}
Ta Thi Hoai An and Nguyen Viet  Phuong, A note on Hayman's conjecture,
 Internat. J. Math. 31 (2020), no. 6, 2050048, 10 pp.
\bibitem{B1}
W. Bergweiler,
On the zeros of certain homogeneous differential polynomials,
Arch. Math. (Basel) 64 (1995), no. 3, 199--202.
\bibitem{BA} W. Bergweiler,
A new proof of the Ahlfors five islands theorem,
J. Anal. Math. 76 (1998), 337--347. 
\bibitem{BN} W. Bergweiler, Normality and exceptional values of derivatives,
Proc. AMS, 129 (2000) 1, 121--129.
\bibitem{BB} W. Bergweiler, Bloch's principle, 
Comput. Methods Funct. Theory 6 (2006), no. 1, 77--108.
\bibitem{BE}
W. Bergweiler and A. Eremenko, On the singularities of the inverse
to a meromorphic function of finite order,
Rev. Mat. Iberoamericana 11 (1995), no. 2, 355--373.
\bibitem{BE2} W. Bergweiler and A. Eremenko, On Bloch's Principle of
Topological Continuity,
\newline
arXiv:2309.14749 to appear in CMFT.
\bibitem{BL}
W. Bergweiler and J. Langley,
Multiplicities in Hayman's alternative,
J. Aust. Math. Soc. 78 (2005), no. 1, 37--57.
\bibitem{BP}
W. Bergweiler and X. Pang, 
On the derivative of meromorphic functions with multiple zeros,
J. Math. Anal. Appl. 278 (2003), no. 2, 285--292.
\bibitem{Bloch} A. Bloch,
La conception actuelle de la th\'eorie des fonctions
enti\`eres et m\'eromorphes, Enseignement Math. 25 (1926), 83--103.
\bibitem{Chang} J. Chang, On meromorphic functions whose first
derivative has finitely many zeros, Bull. LMS, 44 (2012) 703--715.
\bibitem{Chen} H. H. Chen and M. L. Fang,
On the value distribution of $f^nf'$,
Science in China, 38, 7 (1995) 789--798.
\bibitem{C} P. Csillag, \"Uber ganze Funktionen, welche drei nicht
verschwindende Ableitungen besitzen, Math. Ann. 110 (1935), 745--752.
\bibitem{EB} A. Eremenko, Interactions between holomorphic dynamics
and function theory,
\noindent
arXiv:2007.02799.
\bibitem{F} G. Frank and Weissenborn,
Rational deficient functions of meromorphic functions.
Bull. London Math. Soc. 18 (1986), no. 1, 29--33.
\bibitem{H}
W. Hayman, Picard values of meromorphic functions and their derivatives,
 Ann. of Math. (2) 70 (1959), 9--42. 
\bibitem{L1}
J. Langley, Proof of a conjecture of Hayman concerning $f$
and $f^{\prime\prime}$,
J. London Math. Soc. (2) 48 (1993), no. 3, 500--514. 
\bibitem{L2}
J. Langley,
A lower bound for the number of zeros of a meromorphic function
and its second derivative, 
Proc. Edinburgh Math. Soc. (2) 39 (1996), no. 1, 171--185.
\bibitem{L4} J. Langley, The second derivative of a meromorphic
function of finite order, 
Bull. London Math. Soc. 35 (2003), no. 1, 97--108.
\bibitem{Ln}
Wei-Chuan Lin and  Hong-Xun Yi,  On a conjecture of W. Bergweiler,
Proc. Japan Acad. Ser. A Math. Sci. 79 (2003), no. 2, 23–-27
\bibitem{Mc} M. McQuillan, The Bloch Principle, arxiv.org:1209.5402
\bibitem{Mi}
H. Milloux, Sur les fonctions m\'eromorphes et leurs d\'eriv\'ees,
Paris, 1940. 
\bibitem{Miranda}
C. Miranda, Sur un nouveau crit\`ere de normalit\'e pour
les familles de fonctions holomorphes, Bull. Soc. Math. Fr. 63, (1935)
185--196 (1935).
\bibitem{Montel} P. Montel,
Sur les fonctions analytiques qui admettent deux val\'eurs
exceptionnelles dans un domaine, C. R. 153 (1911) 996--998. 
\bibitem{NPZ} S. Nevo, X. Pang, and L.  Zalcman, 
Quasinormality and meromorphic functions with multiple zeros,
Journal d'Analyse, 101 (2007) 1--23.
Electron. Res. Announc. Amer. Math. Soc. 12 (2006), 37–43.
\bibitem{V} L. I. Volkovyskii, Research on the Type Problem of a
Simply Connected Riemann Surface,
Proc. Steklov Math. Inst. 34, Acad. Sci. USSR,
Moscow, 1950.
\bibitem{Pang} Xue Cheng Pang, Bloch’s principle and normal criterion,
Science in China 32
(1989), 782--791.
\bibitem{Picard} E. Picard, Sur les fonctions enti\`eres,
C. R. 89 (1880) 662--665.
\bibitem{Picard1} E. Picard, Sur les fonctions analytiques uniformes
dans le voisinage d’un point singulier essentiels,
C. R. 88 (1879) 745--747. 
\bibitem{Picard0} E. Picard, Sur une propri\'et\'e des fonctions enti\`eres,
C. R. 88 (1879) 1024--1027. 
\bibitem{Schottky} F. Schottky, \"Uber den Picardschen Satz und die Borelschen
Ungleichungen, Sitzungsber. Preuss. Acad. (1904), 1244--1263. 
\bibitem{Schwick} W. Schwick, Normality criteria for families of meromorphic
functions, J. d'Analyse, 52 (1989) 241--289.
\bibitem{Y} K. Yamanoi, 
Zeros of higher derivatives of meromorphic functions in the complex plane,
Proc. Lond. Math. Soc. (3) 106 (2013), no. 4, 703--780.
\bibitem{Yam} K. Yamanoi, Bloch's Principle for holomorphic maps into
subvarieties of semi-Abelian varieties, arXiv:2304.057
\bibitem{Z} L. Zalcman, A heuristic principle in complex function theory, Amer.
Math. Monthly 82 (1975), 813--817.
\end{thebibliography}
\end{document}